\documentclass{amsart}

\usepackage{amsfonts}
\usepackage{cases}
\usepackage{mathrsfs}
\usepackage{bbm}
\usepackage{amssymb}
\usepackage{txfonts}
\usepackage{amscd}
\usepackage{amsfonts,latexsym,amsmath,amsthm,amsxtra,mathdots}
\usepackage[bookmarks,colorlinks]{hyperref}
\usepackage[all,cmtip]{xy}
\RequirePackage{amsmath} \RequirePackage{amssymb}
\usepackage{color}
\usepackage{colordvi}
\usepackage{multicol}
\usepackage{tikz}
\usepackage{hyperref}
\usepackage{graphicx}
\usepackage{amsmath}
\usepackage{amsmath,amscd}
\usetikzlibrary{matrix,arrows}
\usepackage{textcomp}

\makeatletter
\DeclareFontFamily{OMX}{MnSymbolE}{}
\DeclareSymbolFont{MnLargeSymbols}{OMX}{MnSymbolE}{m}{n}
\SetSymbolFont{MnLargeSymbols}{bold}{OMX}{MnSymbolE}{b}{n}
\DeclareFontShape{OMX}{MnSymbolE}{m}{n}{
    <-6>  MnSymbolE5
   <6-7>  MnSymbolE6
   <7-8>  MnSymbolE7
   <8-9>  MnSymbolE8
   <9-10> MnSymbolE9
  <10-12> MnSymbolE10
  <12->   MnSymbolE12
}{}
\DeclareFontShape{OMX}{MnSymbolE}{b}{n}{
    <-6>  MnSymbolE-Bold5
   <6-7>  MnSymbolE-Bold6
   <7-8>  MnSymbolE-Bold7
   <8-9>  MnSymbolE-Bold8
   <9-10> MnSymbolE-Bold9
  <10-12> MnSymbolE-Bold10
  <12->   MnSymbolE-Bold12
}{}
\let\llangle\@undefined
\let\rrangle\@undefined
\DeclareMathDelimiter{\llangle}{\mathopen}%
                     {MnLargeSymbols}{'164}{MnLargeSymbols}{'164}
\DeclareMathDelimiter{\rrangle}{\mathclose}%
                     {MnLargeSymbols}{'171}{MnLargeSymbols}{'171}
\makeatother

\newcommand{\lk}{\operatorname{lk}}

\hypersetup{
    colorlinks,
   linkcolor={red!10!black},
   citecolor={blue!10!black},
    urlcolor={blue!10!black}
}

     \newcommand{\BZ}{{\mathbb {Z}}}

\def\-{^{-1}}

\newcommand{\delete}[1]{}

    \newcommand{\st}{{\mathrm{st}}}

    \theoremstyle{plain}

\newtheorem{thm}{Theorem}[section]
\newtheorem{lem}[thm]{Lemma}
\newtheorem{prop}[thm]{Proposition}
\newtheorem{cor}[thm]{Corollary}
\newtheorem{rem}[thm]{Remark}

\newtheorem*{thmA}{Theorem A}
\newtheorem*{thmB}{Theorem B}
\newtheorem*{thmC}{Theorem C}
\newtheorem*{rem1}{Remark 1}
\newtheorem*{rem2}{Remark 2}

    \numberwithin{equation}{section}

\def\Proof{\noindent{\bf Proof}\quad}
\def\qed{\hfill$\square$\smallskip}

\begin{document}

\title{Poly-freeness of Artin groups and  the Farrell--Jones Conjecture}

\author{Xiaolei Wu}
\address{ Fakult\"at f\"ur Mathematik, Universit\"at Bielefeld, Postfach 100131, D-33501 Bielefeld, Germany}
\email{xwu@math.uni-bielefeld.de}


\date{Jun. 2021}


\begin{abstract}
We provide two simple  proofs of the fact  that even Artin groups of FC-type are  poly-free which was recently established by R. Blasco-Garcia, C. Mart\'inez-P\'erez and L. Paris. More generally, let $\Gamma$ be a finite simplicial graph with all edges labelled by  positive even integers and $A_\Gamma$ be its associated Artin group; our new proof implies  that if $A_T$ is poly-free (resp. normally poly-free) for every clique $T$ in $\Gamma$, then $A_\Gamma$ is poly-free (resp. normally poly-free). We  prove similar results regarding the Farrell--Jones Conjecture for even Artin groups. In particular, we show that if $A_\Gamma$ is an even Artin group such that  each clique in $\Gamma$ either has at most 3 vertices,  has all of its labels  at least $6$, or is the join of these two types of cliques (the  edges connecting the cliques are all labelled by $2$), then $A_\Gamma$ satisfies the Farrell--Jones Conjecture. In addition, our methods enables us to obtain results for general Artin groups. 
\end{abstract}

\maketitle

\section*{Introduction}

 Recently, there has been a growing interest in the study of even Artin groups \cite{BCA18,BMP19,BL19}. In particular, Blasco-Garcia,  Mart\'inez-P\'erez and  Paris proved in \cite{BMP19} that even Artin groups of FC-type are poly-free. One of the goals of this paper is to provide two short proofs of their result. Recall that a group $G$ is called \emph{poly-free} if it admits a chain of subgroups $ 1 =G_0\unlhd G_1\unlhd \cdots \unlhd G_{n-1}\unlhd G_n=G$, such that  $G_i/G_{i-1}$ is a free group of possibly infinite rank. It is called \emph{normally poly-free}  if in addition, each $G_i$ is  normal in $G$.  The minimal such $n$ is called the  poly-free (resp. normally poly-free) length of $G$.  Our first theorem can be stated as follows.

\begin{thmA} \label{main-thm}
Let $\Gamma$ be a finite simplicial graph with all edges labelled by  positive even integers  and $A_\Gamma$ be its associated Artin group. If $A_\Gamma$ is of FC-type, then it is normally poly-free. In general if $A_T$ is poly-free (resp. normally poly-free) for every clique $T$ in $\Gamma$, then $A_\Gamma$ is poly-free (resp. normally poly-free).
\end{thmA}

\begin{rem1}
Blasco-Garcia further proved in \cite{BG20} that any large even Artin group or even Artin group based on a triangle graph is poly-free. Combining this with Theorem A, one obtains  additional Artin groups that are poly-free.
\end{rem1}

In \cite[Question 2]{Be99} and
the discussions below it, Bestvina asks whether all Artin groups are virtually poly-free. In light of Theorem A, it might be true that all even Artin groups are normally poly-free.  Note that our reduction does also work for virtual poly-freeness (see Remark \ref{rem-vpf-eag}).  Poly-free groups have many nice properties, for example, they are torsion-free, locally indicable, have finite asymptotic dimension (one can prove this via induction using \cite[Theorem 2.3]{DS06}) and satisfy the Baum--Connes Conjecture with coefficients \cite[Remark 2]{BKW19}. Furthermore, based on the work of Bestvina, Fujiwara and Wigglesworth \cite{BFW19}, Br\"uck, Kielak and the author showed that normally poly-free groups satisfy the Farrell--Jones Conjecture \cite[Theorem A]{BKW19}. 

Our proof of Theorem A employs Bass-Serre theory. It also implies the following.

\begin{thmB} \label{thm-B}
Let $\Gamma$ be a finite simplicial graph with all edges labelled by  positive even integers  and $A_\Gamma$ be its associated Artin group. If $A_T$ satisfies the Farrell--Jones Conjecture for every clique $T$ in $\Gamma$, then $A_\Gamma$ also does. In particular, if each clique in $\Gamma$ either has at most 3 vertices,  has all of its labels  at least $6$, or is the join of these two types of cliques (the  edges connecting the cliques are all labelled by $2$), then $A_\Gamma$ satisfies the Farrell--Jones Conjecture.
\end{thmB}

Here by Farrell--Jones Conjecture we mean the Farrell--Jones Conjecture with finite wreath products and coefficients in additive categories (which we will abbreviate to FJCw), see \cite{BKW19,Lu18,RV18} and the references therein  for more information. Many Artin groups are known to satisfy FJCw, for example Artin groups of XXL type \cite[Corollary E]{Ha18}, Artin groups of FC-type \cite{HO19} and some classes of affine Artin groups \cite{Ro18, Ro20}.  In the last section, we extend our results to general Artin groups. In particular, we are able to prove the following (Corollary \ref{cor-gel-poly-free}). 

\begin{thmC}
Let $\Gamma_0$ be the join of labelled finite simplicial graphs $\Gamma_1,\cdots, \Gamma_n$ where the connecting edges are labelled by $2$. If for every $i$, either the graph $\Gamma_i$ is a tree or  $A_{\Gamma_i}$ is an even Artin group of FC-type, then for any subgraph $\Gamma\leq \Gamma_0$, $A_\Gamma$ is normally poly-free. In particular, $A_\Gamma$ satisfies FJCw.
\end{thmC}

\begin{rem2}
Theorem C generalizes \cite[Theorem 3.3.5]{Mo19} in Lorenzo Moreschini's Ph.D thesis. Moreover, Lemma \ref{inf-fp-poly-free}  answers his question in \cite[Remark 1.3.12]{Mo19}.
\end{rem2}

\textbf{Acknowledgements.} 
The author was partially supported by Wolfgang L\"uck's ERC Advanced Grant “KL2MG-interactions”
(no. 662400) and the DFG Grant under Germany's Excellence Strategy - GZ 2047/1, Projekt-ID 390685813. He wants to thank Jingyin Huang for some stimulating discussions, and Kai-Uwe Bux, Georges Neaime, Sayed K. Roushon and Shengkui Ye for some helpful comments. He also thanks the anonymous referees for many suggestions on how to improve the
readability of this paper.

\section{Basics about Artin groups}\label{artin-group} 

We give a quick introduction to Artin groups, see \cite{Mc} for more information.  Given a  finite simplicial graph $\Gamma$ with edges $e\in E(\Gamma)$ labelled by integers $m_e \ge 2$, the associated \emph{Artin group} is
	\[
	A_{\Gamma} := \langle V(\Gamma) \mid \langle v,w\rangle^{m_{e}} =\langle w,v\rangle^{m_{e}} \text{ for each } e \in E(\Gamma) \rangle,
	\]
	where $v$ and $w$ are the vertices of $e$, and $\langle v,w \rangle^{m_{e}}$ denotes the alternating product $vwv \cdots$ of length $m_e$ starting with $v$.  For example, $\langle v,w\rangle^3 =vwv$.  An Artin group is called \emph{even} if all edge labels are even.  The group is a right-angled Artin group (RAAG) if all edges are labelled by $2$. Despite the simple presentation, Artin groups in general are poorly understood. For example, the $K(\pi,1)$ Conjecture \cite[Conjecture~1]{CD95} predicts that every Artin group $A_\Gamma$ has a finite classifying space $K(A_\Gamma,1)$ while in general it is even unknown whether $A_\Gamma$ is torsion-free.
	
	The \emph{Coxeter group} $W_\Gamma$ associated to $\Gamma$ is the quotient of $A_\Gamma$ by  the relations $v^2 =1$ for all $v \in V(\Gamma)$.  $A_\Gamma$ is of \emph{spherical type} if the associated Coxeter group $W_\Gamma$ is finite.

	Given a subset $T$ of $V(\Gamma)$, we denote by $A_T$ (resp. $W_T$) the subgroup of $A_{\Gamma}$ (resp. $W_\Gamma$)  generated by $T$, and by $\Gamma_T$ the full subgraph of $\Gamma$ spanned by $T$. The group $W_T$ is  the Coxeter group of $\Gamma_T$ \cite{Bo68}, and $A_T$ is the Artin group of $\Gamma_T$ \cite{Va83}.  We say that $T$ spans a clique in $\Gamma$ if any two elements of $T$ are joined by an edge in $\Gamma$. We say that  $A_\Gamma$ is of FC-type if $A_T$ is of spherical type for each clique in $\Gamma$. Artin groups of FC-type are relatively well understood, for example, the $K(\pi,1)$ Conjecture is known to hold for them \cite{CD95}.

\section{Even Artin groups}
The purpose of this section is to prove our main theorems related to even Artin groups. The proofs use Bass-Serre theory, see  \cite{DD89, Se03} for more information.

We first need the following observation on even Artin groups of FC-type.

\begin{lem}\cite[Lemma 3.1]{BMP19} \label{cha-eafc}
Let $A_\Gamma$ be an even Artin group. Then $A_\Gamma$ is of FC-type if and only if
every triangular subgraph of $\Gamma$ has at least two edges labelled by $2$.
\end{lem}

We also need the following characterization of normally poly-free groups.

\begin{lem}\label{lem-npf-equiv}
A group $G$ is normally poly-free if and only if there exists a sequence of  homomorphisms 
$$ G \xrightarrow[]{q_1}Q_1 \xrightarrow[]{q_2}Q_2 \cdots  \xrightarrow[]{q_n}Q_n=1, $$
such that $\ker(q_i)$ is free for each $i$.
\end{lem}
\Proof 
Suppose $G$ is normally poly-free, i.e. there exists a normal chain 
$$ 1 =G_0\unlhd G_1\unlhd \cdots \unlhd G_{n-1}\unlhd G_n=G,$$ such that  $G_i/G_{i-1}$ is a free group and $G_i$ is normal in $G$. We can choose a sequence of quotient maps as follows.
$$ G =G/G_0\xrightarrow[]{q_1}G/G_1 \xrightarrow[]{q_2}G/G_2 \cdots  \xrightarrow[]{q_n}G/G_n=1.  $$
Then $\ker{q_i} = G_i/G_{i-1}$  is free. For the other direction, up to replace $Q_i$ by $q_i\circ \cdots\circ q_1(G)$ we can assume all the $q_i$s are surjective. Note that $\ker(q_i)$ is still free as subgroups of free groups are again free.  Now let $G_i= \ker(q_1\circ q_2\cdots\circ q_{n-i})$. Then we have $ 1 =G_0\unlhd G_1\unlhd \cdots \unlhd G_{n-1}\unlhd G_n=G$, $G_i\lhd G$ and $G_i/G_{i-1}\cong \ker(q_{n-i})$ which is free.
\qed

Our basic strategy to prove that $A_\Gamma$ is normally poly-free is via mapping $A_\Gamma$ to $A_{\bar\Gamma}$, where $A_{\bar{\Gamma}}$ has the same underlying graph and labels as $\Gamma$ except we change one label to $2$. The following lemma is our key observation. 

\begin{lem} \label{key-lemma}
Let $I_{2n}$ be the connected simplicial graph with one single edge whose label is $2n$ and $A_{I_{2n}}$ be the corresponding Artin group. Then the map 
$$R: A_{I_{2n}} \rightarrow A_{I_2},$$
induced by mapping the vertex generators to the vertex generators is surjective and $\ker(R)$ is a free group.
\end{lem}

\Proof
The map is obviously surjective. The group $A_{I_{2n}}$ has the following presentation:

$$ \langle x,y \mid  (xy)^n = (yx)^n \rangle.$$
The relation $(xy)^n = (yx)^n$ can be rewritten as $(xy)^n = y (xy)^n y^{-1}$. Let $a = xy$ and $t=y$, we have a new presentation of $A_{I_{2n}}$,
$$ \langle a,t \mid  ta^nt^{-1} = a^n \rangle.$$
This is precisely the Baumslag-Solitar group $BS(n,n)$. Now $R$ maps $BS(n,n)$ to the free abelian group of rank $2$ 
$$\langle \bar{x},\bar{y}\mid \bar{x} \bar{y} = \bar{y}\bar{x}  \rangle, $$
via mapping $a$ to $\bar{x}\bar{y}$ and $t$ to $\bar{y}$. Note that $BS(n,n)$ can be viewed as an HNN extension of $\langle a\rangle $ along the subgroup $\langle a^n \rangle $ via the identity map. By Bass-Serre theory,  $BS(n,n)$ acts on a tree $T$ with vertex stabilizers conjugate to $\langle a \rangle$ and edge stabilizers conjugate to $\langle a^n \rangle$ inside $BS(n,n)$. Since $R$ is injective when restricted to $\langle a\rangle$, we have $\ker(R)$ acting on $T$ freely, hence $\ker(R)$ is a free group.
\qed

With Lemma \ref{key-lemma}, we can already show that even Artin groups of spherical type are  normally poly-free. 

\begin{lem}\label{npf-s-type}
Let $\Gamma$ be a clique with all labels even and $A_\Gamma$ be its associated Artin group. If $A_\Gamma$ is of spherical type, then it is a product of copies of $\BZ$ and Artin groups of type $I_{2n}$. In particular, it is normally poly-free.
\end{lem}

\Proof  By the classification of Artin groups of spherical type \cite{Co35}, we know that the only irreducible even Artin groups of finite type are of the form $A_{I_{2n}}$ or $\BZ$. Hence $A_\Gamma$ is a product of copies of $\BZ$ and Artin groups of type $I_{2n}$.  The rest of the lemma now follows from Lemma \ref{key-lemma}.
\qed

\begin{lem}\label{inf-fp-poly-free}
Let $\{G_i\}_{i\in I}$ be a sequence of poly-free groups of length at most $n$. Then the free product $\ast_{i\in I} G_i$ is again a poly-free group of length at most $n$. The analogous statement also holds for normally poly-free groups. Moreover, when $I$ is finite and each $G_i$ is virtually  poly-free (resp. virtually normally poly-free), so is $\ast_{i\in I} G_i$.
\end{lem}
\Proof We first prove the ``poly-free" part of the lemma by induction on $n$. Suppose each $G_i$  is a poly-free group of length at most $n=k+1$, then each $G_i$ has a normal subgroup $H_i$ with poly-free length $k$ and $G_i/H_i$ being a free group. Hence we have a map from $\ast_{i\in I }{G_i}$ to $\ast_{i\in I }G_i/H_i$ which is a free group. Denote the kernel by $K$. By Bass-Serre theory $\ast_{i\in I} G_i$ acts on the corresponding Bass-Serre tree $T$ with vertex stabilizers conjugate  to $G_i$ and edge stabilizers trivial. In particular, $K$ as a subgroup of $\ast_{i\in I }{G_i}$ also acts on $T$. Moreover, since  $K$ is a normal subgroup of $\ast_{i\in I}{G_i}$, it acts on $T$ with vertex stabilizers isomorphic to $H_i$ and edge stabilizers trivial. In fact, given any vertex $v$ in $T$, its stabilizer is of the form $gG_ig^{-1}$ for some $g\in \{G_i\}_{i\in I}$. Restricted to the action of $K$, we have the stabilizer of $v$ is $g G_ig^{-1}\cap K = g(G_i\cap K)g^{-1} = gH_ig^{-1}$ where the first equality uses that $K$ is normal. Thus by Bass-Serre theory, $K$ is a free product of (possibly infinitely many copies of) $H_i$ and some free groups. By induction it is a poly-free group of length  $k$. Hence $\ast_{i\in I} G_i$ is a poly-free group of length $k+1$.

Now when  each $G_i$  is a normally poly-free group of length at most $n=k+1$, we have that each $G_i$ has a free normal subgroup $F_i$ such that $G_i/F_i$ is a normally poly-free group of length at most $k$. We define the map $q:\ast_{i\in I}G_i \to \ast_{i\in I}G_i/F_i$. The same argument implies that $\ker(q)$ is a free group. Now by induction $\ast_{i\in I}G_i/F_i$ is a normally poly-free group of length at most $k$, hence $\ast_{i\in I}G_i$ is a normally poly-free group of length at most $k+1$. 

When $G_i$ is only virtually (normally) poly-free, we can first map $\ast_{i\in I} G_i$  to the direct product $\times_{i\in I} G_i$. Denote the map by $f$, then $\ker{(f)}$ is in fact a free group by the same argument using that $f|_{G_i}$ is injective (see also Remark \ref{rem-graph-grp}). It is  straightforward to check now that $\times_{i\in I} G_i$ is again virtually (normally) poly-free when $I$ is finite. Therefore $\ast_{i\in I} G_i$ is virtually (normally) poly-free.
\qed

\begin{lem}\label{lem-kernel-free} 
Let $C$ be a subgroup of $A$ and $B$ and $G = A\ast_C B$ be the amalgamated product. Suppose $f:G \to Q$ is a map such that $f|_C$ is injective. If   $\ker(f|_A)$  and $\ker(f|_B)$ are free groups (resp. poly-free groups), then $\ker(f)$ is also a free group (resp. poly-free group).
\end{lem}

\Proof We first deal with the case that $\ker(f|_A)$  and $\ker(f|_B)$ are free groups. By Bass-Serre theory, $G$ acts on a tree $T$ such that the vertex stabilizers are either conjugate to $A$ or $B$  and all edge stabilizers are conjugate to $C$ in $G$. Now $\ker(f)$ as a subgroup of $G$ also acts on $T$. For any edge $E\in T$, since its stabilizer is conjugate to $C$ in $G$, we have  $G_E = gCg^{-1}$ for some $g \in G$. This implies that $\ker(f) \cap G_E = (\ker(f) \cap C)^g$ since $\ker(f)$ is a normal subgroup. Now $f|_C$ is injective, we have that $\ker(f) \cap G_E = \{1\}$. For the same reason, we have  $\ker(f) \cap G_v $ is a free group for any vertex $v\in T$. Now this implies that $\ker(f)$ acts on the tree $T$ with vertex stabilizers free and edge stabilizers trivial. By Bass-Serre theory, we have $\ker(f)$ is a free product of free groups and hence it is free. When $\ker(f|_A)$ and $\ker(f|_B)$ are poly-free groups of length at most $d$, the same argument shows that $\ker(f)$ is free product of free groups and poly-free groups of length at most $d$. By Lemma  \ref{inf-fp-poly-free}, we have that $\ker(f)$ is a poly-free group of length at most $d$.
\qed

\begin{rem}\label{rem-graph-grp}
 Lemma \ref{lem-kernel-free} also holds  for graph of groups as long as $f$ restricted to each edge group is injective and the kernel of $f$ restricted to each vertex group is free.
\end{rem}

We are now ready to prove the main theorem. We will first give an almost complete proof using the map in Lemma \ref{key-lemma}.

\textbf{An almost complete proof of Theorem A.} We will show the following: (1), even Artin groups of FC-type are normally poly-free; (2), given an even Artin group $A_\Gamma$, if $A_T$ is poly-free for each clique $T$ in $\Gamma$, then $A_\Gamma$ is poly-free.

Let $A_\Gamma$ be an even Artin group of FC-type. We will prove (1)  by induction on the sum of the number of vertices and the number of edges that are labelled not $2$.  By Lemma \ref{npf-s-type}, we can assume $\Gamma$ is not a clique. In particular, we have a vertex $v_0$ such that the star $\st (v_0) \neq \Gamma$. In this case, let $\Gamma_1$ be the full subgraph of $\Gamma$ spanned by vertices in $V(\Gamma)\setminus{v_0}$, where $V(\Gamma)$ is the vertex set of $\Gamma$. Then we can write $A_\Gamma$ as an amalgamated product of groups, 
$$A_\Gamma =  A_{\st{(v_0)}} \ast_{A_{\lk{(v_0)}}} A_{\Gamma_1}. $$
If all the edges connecting $v_0$ to the link $\lk(v_0)$ are labelled by $2$, we have $A_{\st{(v_0)}} = \BZ \times A_{\lk{(v_0)}}$. We then define a map $f$ from $A_\Gamma$  to $ A_{\Gamma_1}$ by killing $v_0$. Now $f|_{A_{\Gamma_1}}$ is injective and $f|_{A_{\st{(v_0)}}}$ has kernel $\langle v_0\rangle$. Therefore by Lemma \ref{lem-kernel-free}, $\ker(f)$ is a free group. The result now follows from induction.

Assume now there is an edge $I$ connecting $v_0$ which is not labelled by $2$. We define $\bar{\Gamma}$ to be the labelled graph obtained from $\Gamma$ via changing the label of $I$ to $2$ and let $f_R$ be the corresponding map. By induction, $A_{\bar\Gamma}$ is normally poly-free. So it suffices  to show that $\ker{(f_R)}$ is a free group.  Now since $f_R|_{A_{\Gamma_1}}$ is already injective,  by Lemma \ref{lem-kernel-free}, it is enough to show that the kernel of $f_R$ restricted to $A_{\st(v_0)}$ is free. For that let the other end of  $I$  be $w_0$. If $\st(w_0)\cap \st(v_0)\neq \st(v_0)$, then there exists a vertex $w_1$ in $\st(v_0)$ which is not connected to $w_0$ by an edge. Just as before, we can further write $A_{st(v_0)}$ as the following amalgamated product:
$$ A_{\st(v_0)} =  A_{st(w_0)\cap \st(v_0)} \ast_{ A_{\lk(w_0)\cap \st(v_0)}} A_{\Gamma_2},$$
where $\Gamma_2$ is  the full subgraph of $\st(v_0)$ spanned by vertices in $V(\st(v_0))\setminus{w_0}$. Again, $f_R|_{A_{\Gamma_2}}$ is injective, hence  by Lemma \ref{lem-kernel-free}, we only need to show that the kernel of $f_R|_{A_{st(w_0)\cap \st(v_0)}}$ is free.  Now the stars of $v_0$ and $w_0$ in  $\st(v_0)\cap \st(w_0)$ are already the same. Therefore,  we can always assume   that $\st(w_0) \cap \st(v_0) = \st(v_0)$ and let $\Gamma_3$ be the full subgraph of $\st(v_0)$ spanned by vertices other than $v_0$ and $w_0$. Then $\st(v_0)$ is a join of $I$ and $\Gamma_3$. By Lemma \ref{cha-eafc}, all the edges connecting $I$ and $\Gamma_3$ must have label $2$ since $I$ is not labelled by $2$. Thus $A_{\st{(v_0)}} = A_{I} \times A_{\Gamma_3}$ and $\ker(f_R) \cap A_{st(v_0)}$ is a free group by Lemma \ref{key-lemma}. This shows that $\ker(f_R)$ is a free group.

For (2), the proof is even easier as extensions of poly-free groups by poly-free groups are again poly-free. The proof is by induction on the number of vertices. Assume  $\Gamma$ is not a clique, then we can find a vertex $v_0$ such that $\st (v_0) \neq \Gamma$. Just as before we have
$$A_\Gamma =  A_{\st{(v_0)}} \ast_{A_{\lk{(v_0)}}} A_{\Gamma_1}. $$
Since  $A_{\st{(v_0)}}$ is a  subgraph of $\Gamma$, any clique of it is also a clique of $\Gamma$. In particular for any clique $T$ in $A_{\st{(v_0)}}$, $A_T$ is poly-free. By induction, we have $A_{\st{(v_0)}}$ is poly-free. For the same reason $A_{\Gamma_1}$ is poly-free. We define the map $f: A_\Gamma\to A_{\Gamma_1}$ by killing $v_0$ (here we use the fact that all the labels are even). The kernel of $f$ restricted to $A_{\st(v_0)}$ is a subgroup of $A_{\st(v_0)}$, hence it is poly-free. On the other hand, $f$ restricted to $A_{\Gamma_1}$ is injective. By Lemma \ref{lem-kernel-free}, $\ker{(f)}$ is a poly-free group. 
\qed

\begin{rem}\label{rmk-sex-ply}
In general given a short exact sequence of groups $1 \to K\to G\to Q\to 1$ with $K,Q$  normally poly-free, $G$ is not necessarily normally poly-free.  For example, one can take $K = \BZ^2$, $Q=\BZ$ and $G = \BZ^2\rtimes_f \BZ$ with $f = \left( \begin{array}{rr}
2  & 5 \\
1 & 3 \\
\end{array}\right) $. Then $G$ is poly-free but not normally poly-free. This is the reason that our reduction here only works for poly-free groups.
\end{rem}

To get the reduction to also work for normally poly-free groups, we need a slightly different idea.

\begin{prop}\label{prop-key}
Let $\Gamma$ be a finite simplicial graph labelled by even numbers and $A_\Gamma$ be the corresponding Artin group. Then there exists a chain of  maps
$$ A_\Gamma \xrightarrow[]{q_1}A_1  \xrightarrow[]{q_2}A_2 \cdots  \xrightarrow[]{q_n}A_n, $$
and a collection of cliques $T_1,\cdots, T_p$  of $\Gamma$ such that
$\ker{q_i}$ is free for every $i \in \{1,2,\cdots,n\}$, and $A_{n} \cong A_{T_1} \times A_{T_2}\cdots  \times A_{T_p}$.
\end{prop}

\Proof We will make use of the following retraction map that works for any even Artin group, see for example \cite[p.311 Remark (2)]{BMP19}. 

\textbf{Observation.} Given any subset $S$ of the vertex set $V(\Gamma)$,  the inclusion map $A_S\hookrightarrow A_\Gamma$ always admits a retraction $\pi_S: A_{\Gamma} \rightarrow A_S$ which sends $v$ to $v$ if $v\in S$, and sends $w$ to $1$ for  $w\not\in S$. 

Just as in the previous proof, if $\Gamma$ is not a clique, we have a vertex $v_0$ such that $\st (v_0) \neq \Gamma$, and
$$A_\Gamma =  A_{\st{(v_0)}} \ast_{A_{\lk{(v_0)}}} A_{\Gamma_1}, $$
where $\Gamma_1$ is the full subgraph spanned by $V(\Gamma)\setminus v_0$ in $\Gamma$. In particular, we have two retraction maps $\pi_{\st{(v_0)}}: A_\Gamma \rightarrow A_{\st{(v_0)}}$ and $\pi_{\Gamma_1}: A_\Gamma \rightarrow A_{\Gamma_1} $.  Now let 
$$q_1 = \pi_{\st{(v_0)}}\times \pi_{\Gamma_1} : A_\Gamma \rightarrow A_1= A_{\st{(v_0)}} \times A_{\Gamma_1}.$$
Since $\pi_{\st{(v_0)}}$ and $ \pi_{\Gamma_1}$ are retraction maps, we have that the restriction maps of $q_1$ to both ${A_{\st{(v_0)}}} $ and $A_{\Gamma_1} $  are injective. By Lemma \ref{lem-kernel-free}, we have that $\ker(q_1)$ is a free group. Now if both $A_{\st{(v_0)}} $ and $A_{\Gamma_1}$ are cliques,  we are done. Otherwise, we run the same proof to $A_{\st{(v_0)}} $ or $A_{\Gamma_1}$. After finitely many steps, we stop at an Artin group $A_n$ which is a product of Artin groups over cliques in $\Gamma$.
\qed

Theorem A is now an easy corollary of Lemma \ref{lem-npf-equiv} and Proposition \ref{prop-key}.

\textbf{Proof of Theorem A.}  Let $A_\Gamma$ be an even Artin group. By Proposition \ref{prop-key},  there exists a chain of maps
$$ A_\Gamma \xrightarrow[]{q_1}A_1  \xrightarrow[]{q_2}A_2 \cdots  \xrightarrow[]{q_n}A_n, $$
and a collection $T_1,\cdots, T_p$ of cliques of $\Gamma$ such that
$\ker{q_i}$ is free for every $i \in \{1,2,\cdots,n\}$, and $A_{n} \cong A_{T_1} \times A_{T_2}\cdots  \times A_{T_p}$. If $A_\Gamma$ is an even Artin group of FC-type, then by Lemma \ref{npf-s-type}, we have $A_n$ is normally poly-free. So we have further maps 
$$ A_n \xrightarrow[]{q_{n+1}}A_{n+1}  \xrightarrow[]{q_{n+2}}A_{n+2} \cdots  \xrightarrow[]{q_{n+m}}A_{n+m}=1, $$
such that $\ker(q_{n+j})$ is free for every $j$. By Lemma \ref{lem-npf-equiv}, $A_{\Gamma}$ is normally poly-free. This also shows if for each clique $T$ in $\Gamma$, $A_T$ is normally poly-free, then $A_{\Gamma}$ is normally poly-free. The poly-free part is even easier as  extensions of poly-free groups by poly-free groups are again poly-free.
\qed

\begin{rem}\label{rem-vpf-eag}
Note that if all the cliques are only virtually poly-free (resp. virtually normally poly-free), our proof also implies that $A_\Gamma$ is virtually poly-free (resp. virtually normally poly-free). 
\end{rem}

\begin{prop}\label{trg-Artin-FJC}
Let $\Gamma$ be a clique with all edge labels even. If $\Gamma$ has at most $3$ vertices or all labels  are at least $6$, then $A_\Gamma$ is a CAT(0) group. In particular,  $A_\Gamma$ satisfies FJCw.
\end{prop}

\Proof
If $\Gamma$ is a vertex, then $A_\Gamma$ is isomorphic to $\BZ$ and hence CAT(0). If $\Gamma$ is an edge, $A_\Gamma$ is the Baumslag--Solitar group $BS(n,n)$ (by  the proof of Lemma \ref{key-lemma}) which is also known to be CAT(0). If $\Gamma$ is a triangle with labels $p\leq q\leq r$, then $A_\Gamma$ is known to be CAT(0) if $p\geq  4$ \cite[Theorem 4]{BM00}. Now assume $p=2$. If $q=2$, then $A_\Gamma  \cong  \BZ \times A_{I_r}$ which is again CAT(0). Now if $p=2$ and $q\geq 6$, or $p=2,q=4$ and $r\geq 4$, $A_\Gamma$ is CAT(0) by \cite{Ha02}.

If $A_\Gamma$ is a clique with all edges labels at least $5$, it is CAT(0) by \cite{Ha18}. The fact that CAT(0) groups satisfy FJCw is proved in \cite{BL12, We12}. \qed

\textbf{Proof of Theorem B.} Let $A_\Gamma$ be an even Artin group. By Proposition \ref{prop-key},  there exists a chain of  maps
$$ A_\Gamma \xrightarrow[]{q_1}A_1  \xrightarrow[]{q_2}A_2 \cdots  \xrightarrow[]{q_n}A_n, $$
and a collection $T_1,\cdots, T_p$ of cliques of $\Gamma$ such that
$\ker{q_i}$ is free for every $i \in \{1,2,\cdots,n\}$, and $A_{n} \cong A_{T_1} \times A_{T_2}\cdots  \times A_{T_p}$. Since we assumed that $A_{T_i}$ satisfies the FJCw, the direct product will also satisfy FJCw \cite[Theorem 1.1 (2)]{BKW19}. Since any free-by-cyclic group satisfies FJCw \cite{BKW19}, the first part of Theorem B now follows from  \cite[Theorem 1.1 (4)]{BKW19} and induction.

The second part of the theorem  follows from Proposition \ref{trg-Artin-FJC}. Note that if $T$ is the join of $T_1$ and $T_2$ such that all the edges connecting $T_1$ and $T_2$ are labelled by $2$, then $A_T$ is the direct product of $A_{T_1}$ and $A_{T_2}$. 
\qed

\section{ Artin groups in general}
We give some further results concerning the  poly-freeness of general Artin groups in this section.

\begin{lem}\label{pf-add-edge}
Let $\Gamma_1 \subseteq \Gamma_2$ be  labelled finite  simplicial graphs such that $V(\Gamma_1) =V(\Gamma_2)$, $\Gamma_2$ is obtained from $\Gamma_1$ by adding an edge, and $f: A_{\Gamma_1} \to A_{\Gamma_2}$ be the induced map on Artin groups. Then $\ker(f)$ is a free group.
\end{lem}
\Proof 
Let  $v,w$ be the vertices of the extra edge $e$ and let $\Gamma_0$ be the subgraph of $\Gamma_1$ spanned by the vertices in $V_{\Gamma_1}\setminus{v}$, where $V_{\Gamma_1}$ is the vertex set of $\Gamma_1$.  Then we can write $A_{\Gamma_1}$ as an amalgamated product of groups 
$$A_{\Gamma_1} =  A_{\st{(v)}} \ast_{A_{\lk{(v)}}} A_{\Gamma_0}. $$
Since $\st{(v)}$ and $\Gamma_0$ are full subgraphs of $\Gamma_1$, we have $f|_{ A_{\st{(v)}}}$ and $f|_{A_{\Gamma_0}}  $ are injective. Hence by Lemma \ref{lem-kernel-free}, $\ker{(f)}$ is a free group. \qed

\begin{lem}\label{lem-subg-plf}
Let $\Gamma_1 \subseteq \Gamma_2$ be  labelled finite  simplicial graphs and $A_{\Gamma_1}, A_{\Gamma_2}$ be the corresponding Artin groups. If $ A_{\Gamma_2}$ is a (virtually) poly-free group (resp. normally poly-free), so is  $A_{\Gamma_1}$.
\end{lem}
\Proof Since poly-freeness (resp. normally poly-freeness) passes to subgroups, we can assume  $\Gamma_2$ has the same number of vertices of $\Gamma_1$. Indeed we can take $\Gamma_2'$ to be the full subgraph of $\Gamma_2$ spanned by the vertices of $\Gamma_1$ and let $A_{\Gamma_2'}$ be the corresponding Artin group. Then $A_{\Gamma_2'}$ is a subgroup of $ A_{\Gamma_2}$. The lemma now follows from Lemma \ref{pf-add-edge} by  induction.  \qed

\begin{thm}\begin{enumerate}\label{Cor-AG}
\item If  every Artin group over a clique is (virtually) poly-free  (resp. normally poly-free), then all Artin groups are (virtually) poly-free (resp. normally poly-free).
    \item If every Artin group over a clique is torsion-free, then all Artin groups are torsion-free.
    \item \label{Cor-AG-FJC-1} Let $\Gamma_1 \subseteq \Gamma_2$ be  labelled finite  simplicial graphs and $A_{\Gamma_1}, A_{\Gamma_2}$ be the corresponding Artin groups. If $ A_{\Gamma_2}$ satisfies FJCw, so does  $A_{\Gamma_1}$.
    
    \item \label{Cor-AG-FJC-2} If every Artin group over a clique satisfies FJCw, then all Artin groups satisfy FJCw.
\end{enumerate}

\end{thm}

\Proof 
Item (a) and (b) follow easily from Lemma \ref{pf-add-edge} and Lemma \ref{lem-subg-plf}. We indicate how to prove (c) and (d). Since FJCw is closed under taken subgroups \cite[Theorem 1.1(1)]{BKW19}, we can assume $V(\Gamma_1) =V(\Gamma_2)$ by taking $\Gamma_2$ to be the full subgraph of bigger graph spanned by the vertices in $\Gamma_1$. If $\Gamma_1 = \Gamma_2$, we are done. Otherwise we first assume that $\Gamma_1$ is obtained from $\Gamma_2$ by adding one single edge. Let $q:\Gamma_1\to \Gamma_2$ be the corresponding map. By Lemma \ref{pf-add-edge}, $\ker(f)$ is a free group which satisfies FJCw. By \cite[Theorem 1.1(4)]{BKW19}, we only need to show for each infinite cyclic subgroup $C$ of $A_{\Gamma_2}$, $f^{-1}(C)$ satisfies FJCw. This follows from \cite[Theorem A]{BKW19} since $f^{-1}(C)$ is a free-by-cyclic group. The general case follows from induction.

Now given any Artin group $A_\Gamma$ which is not a clique, we can add edges to $\Gamma$ with label $2$ to make it a clique which we denote by $T$. By (c), if we know FJCw holds for $A_T$, it also holds for $A_\Gamma$. 
\qed

\begin{lem} \label{key-lemma-2}
Let $I_{n}$ be an edge with label $n\geq 2$ and $A_{I_{n}}$ be the corresponding Artin group. Then the map 
$$\chi: A_{I_{n}} \rightarrow \BZ,$$
induced by mapping all the vertex generators to $1\in \BZ$ is surjective and $\ker(\chi)$ is a free group. In particular $A_{I_{n}}$ is normally poly-free.
\end{lem}
\Proof 
When $n$ is even, the proof follows that of  Lemma \ref{key-lemma}. We explain the proof for $n$ odd\footnote{The author learned this part of the proof from Shengkui Ye.}. Recall $A_{I_{2k+1}} = \langle x,y\mid x(yx)^k = (yx)^ky \rangle$. Let $s= x(yx)^k$ and $t = yx$. Then $x = st^{-k}$ and $y = t^{k+1}s^{-1}$. Therefore, the relator $x(yx)^k = (yx)^ky$ becomes 

$$s= t^k t^{k+1}s^{-1},$$ 
which is $s^2 = t^{2k+1}$. Thus we have the following presentation for $A_{I_{2k+1}}$,
$$ \langle s,t\mid  s^2 = t^{2k+1} \rangle. $$
Hence $A_{I_{2k+1}}$ can be viewed as the amalgamated product of $\langle s\rangle$ and $\langle t\rangle$ identifying the subgroup $\langle s^2\rangle$ with $\langle t^{2k+1}\rangle$. 
Now $\chi(s) = 2k+1$ and $\chi(t) = 2$. We have that both  $\chi|_{\langle s\rangle}$ and $\chi|_{\langle t\rangle}$ are injective. Hence by Lemma \ref{lem-kernel-free}, $\ker(\chi)$ is a free group. 
\qed

\begin{rem}
Alternatively, note that  $A_{I_{n}}$ is actually  the fundamental group of a $(2, n)$-torus link complement (see for example \cite[Introduction]{MMV01}). By \cite[Corollary A.1]{MMV01}, $\ker{(\chi)}$ is finitely generated. By Stallings fibration Theorem \cite{St62},  the $3$-manifold fibers over a circle with fiber some surface with boundary and $\chi$ is the induced map on $\pi_1$ for the bundle projection map. This in particular implies that  $\ker{(\chi)}$ is in fact a finitely generated free group. 
\end{rem}

\begin{prop}\label{npf-tree}
Let $\Gamma$ be a labelled finite simplicial graph and $A_{\Gamma}$ be the corresponding Artin group. If $\Gamma$ is a tree, then $A_{\Gamma}$ is normally poly-free.
\end{prop}
\Proof We prove the proposition by induction. If $\Gamma$ has only one edge, $A_{\Gamma}$ is normally poly-free by Lemma \ref{key-lemma-2}. Suppose $\Gamma$ has more than one edge, and let $v_0$ be a vertex of valence one and $e$ be the unique edge in $\Gamma$ connecting $v_0$. Let $\Gamma_1$ be the full subgraph of $\Gamma$ with vertex set $V(\Gamma)\setminus v_0$. Then we have 
$$A_\Gamma = A_{e}  \ast_{A_{w_0}} A_{\Gamma_1}  $$
where $w_0$ is the other vertex of $e$. Now we define a map $q:A_\Gamma \to A_{\Gamma_1}$ as follows: $q(v)=v$ if $v\neq v_0$ and $q(v_0)=w_0$. In particular, $q|_{A_{\Gamma_1}}$ is injective and the kernel of $q|_{A_e}$ is a free group. Hence, by Lemma \ref{lem-kernel-free}, $\ker(q)$ is a free group. By induction, $A_{\Gamma_1}$ is normally poly-free. So $A_\Gamma$ is normally poly-free.
\qed

Combining Theorem A, Lemma \ref{lem-subg-plf}, Proposition \ref{npf-tree} and \cite[Theorem A]{BKW19}, we have the following. 
\begin{cor}\label{cor-gel-poly-free}
Let $\Gamma_0$ be the join of labelled finite simplicial graphs $\Gamma_1,\cdots, \Gamma_n$ where the connecting edges are labelled by $2$. If for every $i$, either the graph $\Gamma_i$ is a tree or  $A_{\Gamma_i}$ is an even Artin group of FC-type, then for any subgraph $\Gamma\leq \Gamma_0$, $A_\Gamma$ is normally poly-free. In particular, $A_\Gamma$ satisfies FJCw.
\end{cor}

\end{document}